\documentclass{article}
\usepackage{amsfonts}
\usepackage{amsmath}
\usepackage[latin1]{inputenc}
\usepackage[T1]{fontenc}
\usepackage{lmodern}
\usepackage[frenchb]{babel}
\usepackage{graphicx}

\newtheorem{lemme}{Lemme}
\newtheorem{theoreme}{Théorème}
\newtheorem{prop}{Proposition}
\newtheorem{corollaire}{Corollaire}

\title{Les groupes de triangles $(2,p,q)$ sont déterminés par leur spectre des longueurs}
\author{\small{Emmanuel Philippe\footnote{e-mail: matahari@netcourrier.com}}\\\small{Laboratoire Emile Picard}\\\small{Université Paul Sabatier}\\\small{Toulouse}}
\date{}
\begin{document}
\maketitle
\begin{abstract}
On décrit le début du spectre des longueurs des groupes de triangles ayant un angle droit et on montre que le spectre des longueurs caractérise la classe d'isométrie d'un tel groupe.\\
$2000$ $Mathematics$ $subject$ $classification$: 20H10, 32G15, 53C22
\end{abstract}
\begin{center}
\subsubsection*{Abstract}
We describe the length spectra of triangle groups $(2,p,q)$ before
showing that the length spectra characterizes the isometry class of
such a group.\\
$2000$ $Mathematics$ $subject$ $classification$: 20H10, 32G15,
53C22.
\end{center}

\tableofcontents

\vspace{1cm} Si $S$ est une surface hyperbolique (i.e. de
caractéristique d'Euler strictement négative) fermée (i.e.
compacte sans bord), on peut s'intéresser à son spectre des
longueurs, c'est à dire à l'ensemble des longueurs des géodésiques
fermées de $S$ rangées dans l'ordre croissant en respectant les
multiplicités eventuelles. Cet ensemble considéré, on peut se
demander l'information géométrique qu'il contient: peut-on déduire
du
spectre des longueurs le genre de $S$? la classe d'isométrie de $S$?\\
Dans le cas où $S$ est fermée, ces questions on été largement
traitées. En particulier, on montre que la donnée du spectre des
longueurs et équivalente à la donnée du spectre du Laplacien sur
$S$ (c'est le théorème de Huber, \cite{Buser1}). On peut déduire
de cela que le spectre des longueurs détermine l'aire de la
surface, et donc le genre ou, ce qui revient au même, la classe
d'homéomorphie (on pourra se référer à \cite{BGM}). Pour ce qui
est de la classe d'isométrie de la surface, la réponse est
négative, même si l'on dispose du théorème de
Wolpert(\cite{Wolpert},\cite{Buser1}) qui informe qu'une surface
hyperbolique fermée de genre $g$ est \og en général \fg (il
convient de donner un sens précis à cette formule\dots) déterminée
à isométrie près par son spectre des longueurs. On sait notamment
construire explicitement des surfaces de genre $g\geq 4$ ayant le
même spectre des longueurs et étant géométriquement distinctes (on
renvoie à \cite{Buser1} pour une présentation
élégante et synthétique de ces exemples).\\
On dispose, parallélement, de certains résultats de rigidité, par
exemple pour les surfaces hyperboliques compactes de genre $1$
avec une seule composante de bord (\cite{Buser2}), pour les
pantalons compacts ou avec 1,2 ou 3 cusps (\cite{Haas}).

Dans le cas où $S$ est fermée et de genre 2 ou 3, la question reste ouverte.\\

Si maintenant $S$ admet des points coniques (i.e. si son groupe
fondamental contient des éléments elliptiques), on définie de la
même manière le spectre des longueurs, en demandant aux
géodésiques d'être les courbes qui minimisent la distance entre
leurs points sur le complémentaire des points coniques. On dispose
d'un analogue partiel du théorème de Huber dans ce cas (on
consultera \cite{Dryden}).

On se propose ici de démontrer que pour les surfaces hyperboliques
de genre $0$ ayant trois points coniques dont l'un est d'ordre
$2$, le spectre des longueurs permet de déterminer l'ordre de tous
les points coniques, et donc la classe d'isométrie. Pour des
calculs numériques de sytoles dans certains groupes de triangles,
on renvoie à \cite{HK},\cite{LW} et \cite{Vogeler}, cette dernière
référence s'attardant sur le groupe $\Gamma(2,3,7)$.

Ceci est un premier résultat de rigidité pour le cas des surfaces
hyperboliques fermées à points coniques, qui s'applique en
particulier à la surface modulaire.

\newpage

\section{Préliminaires}
\subsection{Introduction}
On considère ici le demi-plan de Poincaré
\[\mathbf{H}=\{ z\in\mathbf{C}~;~\mathrm{Im} z >0 \} \]
que l'on munit de la distance hyperbolique
\[\cosh d(x,y)=1+\frac{|x-y|^2}{2~\mathrm{Im} x~\mathrm{Im} y } \]
Cet espace possède un groupe d'isométrie
\textrm{Isom}~$\mathbf{H}$ qui s'identifie à $PGL(2,\mathbf{R})$
tandis que le groupe des isométries directes est assimilé à
$PSL(2,\mathbf{R})$. Un élément $\gamma$ de
\textrm{Isom}$^{+}\mathbf{H}$ est alors hyperbolique si la trace
d'une des matrices associées de $SL(2,\mathbf{R})$ est strictement
plus grande que $2$ en valeur absolue. Un tel élément agit sur
$\mathbf{H}$ en laissant stable une géodésique appelée axe de
$\gamma$ et le long de laquelle $\gamma$ agit comme une
translation de distance de translation
\[ l(\gamma)=\mathrm{Inf}\{d(x,\gamma x)~;~x\in\mathbf{H}\} \]

Si maintenant $2\leq p\leq q$ sont des entiers vérifiant
\[\frac1p+\frac1q<\frac12 \]
on leur associe un triangle $T$ de $\mathbf{H}$ ayant pour angles
$\pi/2,\pi/p,\pi/q$. Le groupe d'isométrie engendré par les
reflexions relatives aux côtés de $T$ est noté $\Gamma_0(2,p,q)$
tandis que l'ensemble de ses éléments préservant l'orientation est
appelé $\Gamma(2,p,q)$. C'est ce dernier groupe que l'on va
étudier ici.

$\Gamma_0(2,p,q)$ agit sur $\mathbf{H}$ avec $T$ comme domaine
fondamental, ce qui fournit un pavage $\mathcal{P}_0$ admettant
trois types de sommets: ceux de valence $2$, de valence $p$ ou de
valence $q$. Considérons $\mathcal{P}$ le sous pavage obtenu en
considérant uniquement les sommets de valence $q$ et les arètes
les reliant: celui-ci est constitué de $p$-gones de côté $2c$ avec

\[\cosh c =\frac{\cos\frac{\pi}p}{\sin\frac{\pi}q}
\]

et les angles aux sommets sont tous égaux à $2\pi/q$.

Nous dirons que $\Gamma(2,p,q)$ et $\Gamma(2,p',q')$ sont
isométriques si les quotients $\mathbf{H}\backslash\Gamma(2,p,q)$
et $\mathbf{H}\backslash\Gamma(2,p',q')$ le sont.\\

On décrit ici le début du spectre des longueurs
\[ \mathrm{Lsp}~\Gamma(2,p,q)=\{l(\gamma)~;~\gamma \textrm{ élément
hyperbolique de } \Gamma(2,p,q) \} \] qui est l'ensemble des
distances de translation comptées avec multiplicités et ordonnées
dans l'ordre croissant. Le plus petit élément de cet ensemble est
la systole.\\

Commençons par expliquer notre approche, avant de donner la liste
des résultats prouvés dans cet article.\\

Il est naturel de penser qu'un élément hyperbolique $\gamma$ de
$\Gamma(2,p,q)$ ayant une distance de translation petite (disons
$\leq l_0$) déplace au moins un sommet de $\mathcal{P}$ d'un
distance petite également (Idée 1).

De plus, on peut penser que deux points du pavage proches pour la
distance hyperbolique le restent également si on considère leur
éloignement au sens de la distance combinatoire $D$ associée au
pavage (cette distance étant le nombre minimal d'arètes permettant
de relier ces deux points) (Idée 2).

Dans la pratique, nous aurons donc besoin

\begin{itemize} \item de
définir le \og déplacement minimal \fg de $\gamma$ sur
$\mathcal{P}$, noté $\lambda(\gamma)$; \item de relier la distance
hyperbolique $d$ et la distance combinatoire  $D$
\end{itemize}

Le résultat obtenu (propriété \ref{p_1}) s'écrit alors
\[\forall l_0,\exists n_0~;~\forall \gamma \textrm{ hyperbolique },~l(\gamma)\leq l_0
\Rightarrow \lambda(\gamma)\leq n_0 \label{eq1} \]

La démonstration de ce fait repose sur deux résultats formalisant
les idées précédentes (lemme \ref{l_1} et lemme \ref{l_2}).

Le premier indique qu'effectivement le \og rayon hyperbolique \fg
$\rho$ des boules combinatoire augmente avec le rayon combinatoire
de celles-ci.

Le second affirme comme il était pressenti que tout élément
hyperbolique $\gamma$ de distance de translation $l(\gamma)$
majorée par $l_0$ déplacera un point du pavage $\mathcal{P}$ d'une
distance majorée par
\[ C(l_0)=\mathrm{Argcosh}~[(\cosh c)^2(\cosh l_0 -1)+1] \]
où cette constante ne dépend que de $l_0$ et des entiers $p,q$.

On obtient alors la propriété \ref{p_1} cherchée.\\

Nous utilisons ensuite les résultats de la section \ref{s_2} pour
\begin{itemize}
\item décrire le début du spectre des longueurs de
$\Gamma(2,p,q)$; \item en déduire que le spectre des longueurs
caractérise le groupe $\Gamma(2,p,q)$ parmi les groupes de
triangles de ce type.
\end{itemize}

On a le résultat de rigidité spectral suivant (théorème \ref{t_1}):\\

\textbf{Théorème~A:}\\
$Soit$ $\Gamma(2,p,q)$ $et$ $\Gamma(2,p',q')$ $deux$ $groupes$
$de$ $triangles$ $isospectraux$ $au$ $sens$ $des$ $longueurs.$
$Alors$ $ces$ $deux$ $groupes$ $sont$ $isom\acute{e}triques.$\\

Au passage, on aura démontré et exploité le résultat suivant
(corollaire \ref{c_2}, section \ref{s_51}):\\

\textbf{Théorème~B:}\\
$La$ $systole$ $de$ $\Gamma(2,p,q)$ $est$
$d\acute{e}termin\acute{e}e$ $par$ $l'alternative$ $suivante:$
\begin{enumerate}
\item $Si$ $p\geq 4$, $il$ $s'agit$ $de$
$2~\mathrm{Argcosh}~[2\cos\frac{\pi}p\cos\frac{\pi}q]$ \item $Si$
$p=3$, $il$ $s'agit$ $de$
$2~\mathrm{Argcosh}~[2(\cos\frac{\pi}q)^2-\frac12]$
\end{enumerate}

La section \ref{s_2} expose comme annoncé la preuve de la
propriété \ref{p_1}. C'est l'aspect théorique de l'étude.

La section \ref{s_3} est consacrée au calcul explicite des
distances de translation des hyperboliques vérifiant
$\lambda(\gamma)=1,2$.

La section \ref{s_4} est dédiée au calcul de $\rho(3)$ et à la
mise en pratique de la propriété \ref{p_1} pour les longueurs qui
nous intéressent.

La section \ref{s_5} démontre le résultat de rigidité spectral.

\subsection{Notations}
Rappelons que tous les points du pavage $\mathcal{P}$ sont par
hypothèse de valence $q$.

Si $x,y,z$ sont trois éléments de $\mathbf{H}$, on désigne par
$\angle(xy,xz)$ l'angle orienté formé par les segments géodésiques
$xy$ et $xz$ en $x$. En l'abscence de précisions, on considère la
mesure de cet angle dans $[0,2\pi[$. La mesure géométrique de cet
angle est la valeur absolue de sa mesure principale.

Si $x,y$ sont deux sommets de $\mathcal{P}$, on considère
l'ensemble des chemins $\beta$ du pavage formés d'arètes
consécutives de $\mathcal{P}$ et reliant $x$ et $y$. Le nombre
d'arètes d'un tel chemin est sa longueur, notée $L(\beta)$. La
distance combinatoire entre $x$ et $y$ est
\[ D(x,y)=\textrm{Inf}\{ L(\beta)~;~\beta\textrm{ reliant $x$ et $y$
}  \} \] C'est un entier.

Fixons maintenant $s_0$ dans $\mathcal{P}$. Pour $n\geq 1$, la
boule combinatoire de centre $s_0$ et de rayon combinatoire $n$
est l'ensemble des sommets de $\mathcal{P}$ qui se trouvent à une
distance combinatoire au plus égale à $n$ de $s_0$. On pose alors
\[ \rho(n)=\mathrm{Inf}\{d(s_0,s)~;~D(s_0,s)=n\} \]
Moralement, il s'agit du rayon hyperbolique minimal de la sphère
combinatoire de rayon $n$. Ce nombre est évidemment indépendant de
$s_0$. Nous sommes notamment amenés, pour démontrer les résultats
annoncés, à calculer $\rho(3)$ dans la section \ref{s_41}.

Décrivons de manière plus minutieuse la sphère combinatoire de
rayon $3$. Si $s$ vérifie $D(s_0,s)=3$, il existe un chemin
$\beta=\{s_0s_1,s_1s_2,s_2s\}$ du pavage reliant $s_0$ à $s$. On
note $2k_1\pi/q$ et $2k_2\pi/q$ les mesures dans $[0,2\pi[$ de
l'angle orienté de $\beta$ en $s_1$ et $s_2$. Le sommet $s$
s'écrit alors $s(k_1,k_2)$ et on note
\[ d(k_1,k_2):=d(s_0,s(k_1,k_2))\]
Remarquons qu'il existe plusieurs points du type $s(k_1,k_2)$ pour
des entiers fixés et qu'ils sont tous à la même distance
hyperbolique de $s_0$. C'est cette distance, invariante, que l'on
étudie: ceci justifie la notation précédente.

Enfin, pour tout $\gamma$ élément hyperbolique de $\Gamma(2,p,q)$,
on définit son niveau
\[ \lambda(\gamma)=\textrm{Inf}\{~D(s,\gamma s)~;~s\in\mathcal{P}\} \]
qui mesure le déplacement minimal de $\gamma$ sur $\mathcal{P}$ au
sens combinatoire du terme.

\section{Aspect théorique de l'étude}\label{s_2}

Il s'agit ici de démontrer la propriété \ref{p_1} annoncée en
introduction, sur laquelle repose l'étude pratique du spectre des
longueurs.

Començons par montrer le lemme suivant établissant la croissance
de $\rho:\mathbf{N}^*\longrightarrow\mathbf{R}^*$ et reliant ainsi
les distances $d$ et $D$.

\begin{lemme}\label{l_1}
L'application $\rho$ est croissante
\end{lemme}

\textbf{Preuve du lemme \ref{l_1}:}\\
Fixons $s_0$ un sommet de $\mathcal{P}$ et $n\geq2$ un entier.
Soit $x$ un sommet  à distance combinatoire $n$ de $s_0$. Montrons
qu'il existe toujours un élément $y$ de $\mathcal{P}$ étant à une
distance combinatoire $n-1$ de $s_0$ et qui soit plus proche de
$s_0$ que $x$ au sens de la distance hyperbolique. En choisissant
$x=x_0$ de manière adaptée, on obtiendra alors $y=y_0$ vérifiant
\[ \rho(n)=d(s_0,x_0)\geq d(s_0,y_0)\geq\rho(n-1) \] et le lemme sera
démontré.

Soit donc $x$ avec la propriété invoquée. Il existe un chemin
$\beta=\{c_1,\dots,c_n\}$ du pavage reliant $s_0$ à $x$. On écrit
$c_n=xy$ et on désigne par $z$ le milieu de $xy$, qui est un point
de valence $2$ du pavage $\mathcal{P}_0$. L'élément $y$ est de
valence $q$ et vérifie
\[ d(s_0,x)\geq d(s_0,y) \]
En effet, si $m$ est la médiatrice de $xy$, $m$ sépare les points
$x$ et $s_0$: établissons ce fait.

Commençons par écrire le complémentaire de $m$ comme union de deux
sous espaces ouverts $H_x$ et $H_y$, avec des notations évidentes,
et supposons que $s_0$ est dans $H_x$. Nous allons contredire la
minimalité de $\beta$.

Le chemin $\beta$ intersecte $m$ en au moins deux points, dont
l'un est $z$. Notons $w$ le premier point d'intersection de
$\beta$ avec $m$ quand $\beta$ est parcouru de $s_0$ vers $x$.
C'est un sommet de $\mathcal{P}_0$, donc de valence $q$ ou $2$.\\

\underline{\textit{Cas 1: $w$ est de valence $q$}}\\
On écrit tout d'abord \[ \beta=\beta_1\cup\beta_2\cup\{yx\} \] où
$\beta_1$ relie $s_0$ à $w$ et $\beta_2$ relie $w$ à $y$. On pose
alors
\[ \beta'=\beta_1\cup\overline{\beta_2} \] où $\overline{\beta_2}$
désigne le chemin de $\mathcal{P}$ image de $\beta_2$ par la
reflexion d'axe $m$. Le chemin $\beta'$ relie $s_0$ à $x$ et se
trouve être de longueur strictement plus petite que $\beta$, ce
qui contredit la minimalité de $\beta$ et impose à $s_0$
d'appartenir à $H_y$.\\

\underline{\textit{Cas 2: $w$ est de valence $2$}}\\
$w$ est alors le milieu d'une arête $x'y'$ avec $x'$ et $y'$ dans
$\mathcal{P}$. Nous supposons que $x'\in H_x$ et écrivons
\[ \beta=\beta_1\cup\beta_2\cup\{yx\}\cup\{x'y'\} \]
où $\beta_1$ relie $s_0$ à $x'$ et $\beta_2$ relie $y'$ à $y$. Il
suffit maintenant de considérer
\[ \beta'=\beta_1\cup\overline{\beta_2} \]
pour contredire là encore la minimalité de $\beta$.\\

Dans tous les cas, on a  bien $s_0\in H_y$ et
\[ d(s_0,x)\geq d(s_0,y) \]

Reste à constater que $y$ est exactement à une distance
combinatoire $n-1$ de $s_0$: en effet,
$\beta'=\{c_1,\dots,c_{n-1}\}$ est un chemin reliant $s_0$ à $y$
(ce qui fournit une majoration de la distance cherchée) et si
$D(s_0,x)$ était strictement inférieure à $n-1$, on aurait $\beta$
de longueur strictement inférieure à $n$.

Ceci achève la preuve du lemme \ref{l_1}. $\square$

Montrons maintenant le lemme \ref{l_2} qui assure l'existence d'un
\og déplacement minimal \fg maximal pour les éléments
hyperboliques dont la distance de translation est bornée.

\begin{lemme}\label{l_2}
Soit $l_0>0$. Alors il existe une constante $C(l_0)$ telle que
tout élément hyperbolique $\gamma$ de $\Gamma(2,p,q)$ ayant une
distance de translation inférieure ou égale à $l_0$ déplace au
moins un sommet $s_0$ de $\mathcal{P}$ d'une distance
\[ d(s_0,\gamma s_0)\leq C(l_0) \]On a de plus
\[ \cosh C(l_0)= (\cosh c)^2(\cosh l_0 -1) +1 \]
\end{lemme}

\textbf{Preuve du lemme \ref{l_2}:}\\
Considérons $ss'$ une arète du pavage $\mathcal{P}$ qui possède
une intersection non vide avec $A$, l'axe de $\gamma$. Si $s$ est
l'extrémité de l'arète la plus proche de $A$ au sens hyperbolique,
cette distance est nécessairement inférieure à $c$. On obtient
alors
\begin{align*}
\cosh d(s,\gamma s) &=(\cosh d(s,A))^2(\cos l(\gamma)-1)+1 \\&\leq
(\cosh c)^2(\cosh l(\gamma)-1)+1
\\ &\leq (\cosh c)^2(\cosh l_0 -1)+1
\end{align*}
On peut donc prendre $s_0=s$.$\square$\\

On est maintenant en mesure de montrer la
\begin{prop}\label{p_1}
Soit $l_0>0$. Alors il existe un entier $n_0$ tel que, pour tout
élément hyperbolique $\gamma$ de $\Gamma(2,p,q)$,
\[ l(\gamma)\leq l_0 \Rightarrow \lambda(\gamma) \leq n_0 \]
\end{prop}

\textbf{Preuve de la proposition \ref{p_1}:}\\
Par croissance de $\rho$, il existe un plus petit entier $n$ tel
que $\rho(n)>C(l_0)$. Supposons dès lors que $\lambda(\gamma)>n$.
Pour tout $s$ dans $\mathcal{P}$, on a
\[ d(s,\gamma s)\geq\rho(\lambda(\gamma))\geq\rho(n)>C(l_0) \]
et le lemme \ref{l_1} permet alors d'affirmer que $l(\gamma)>l_0$.
Nous pouvons donc choisir $n_0=n-1$. $\square$

\section{Aspect pratique de l'étude: première
partie}\label{s_3}

On décrit ici les valeurs du spectre des longueurs apparaissant
avec les éléments hyperboliques de niveau $1$ et $2$, c'est à dire
les distances de translations correspondant aux éléments
hyperboliques de $\Gamma(2,p,q)$ qui déplacent au moins un point
du pavage d'une distance combinatoire $1$ ou $2$. On y trouve des
formules explicites élémentaires et on ordonne ces valeurs dans
l'ordre croissant, en tenant compte des valeurs des entiers $p$ et
$q$.

\subsection{Les hyperboliques de niveau 1} \label{s_31}

Soit $\gamma$ un élément hyperbolique de $\Gamma(2,p,q)$ déplacant
au moins un sommet $s_0$ d'une distance combinatoire $1$.
Choisissons un $s_0$ tel que $D(s_0,\gamma s_0)=1$ et décomposons
$\gamma$ en produit de deux réflexions dont les axes sont des
géodésiques $G$ et $G'$ du
pavage $\mathcal{P}_{0}$: expliquons comment.\\

Si $s_0$  et tel que $D(s_0,\gamma s_0)=1$, écrivons
$\gamma^2s_0=x$. L'arète $\gamma s_0x$ se trouve donc être l'image
par $\gamma$ de l'arète $s_0\gamma s_0$. Prenons pour $G$ la
médiatrice de $s_0\gamma s_0$ et pour $G'$ la bissectrice de
l'angle formé par ces deux arètes au sommet $\gamma s_0$. On note
$G'=G'_k$ où $k$ est un entier entre $1$ et $q-1$ tel que
$\angle(\gamma s_0 s_0,\gamma s_0 x)=2k\pi/q$.  Un calcul
trigonométrique permet d'affirmer que la distance de translation
de $\gamma=\bar{r}_{G'_k}\bar{r}_G$ est alors donnée par
\[l_1(k)= 2~\mathrm{Argcosh}~D_1(k)~;~ 1 \leq k \leq q-1 \]
dès que la quantité
\[D_1(k)=\sin\frac{k\pi}q \frac{\cos\frac{\pi}p}{\sin\frac{\pi}q}~;~1 \leq k \leq
q-1\] est strictement plus grande que 1: dans le cas contraire,
$\gamma$ n'est pas hyperbolique (\cite{Beardon}, p. 157 et p.
179).

On constate que les valeurs prises par $D_1(k)$ sont symétriques
par rapport à la valeur en la partie entière de $q/2$. On a
notamment:
\begin{align*}
D_1(1)&=\cos\frac{\pi}q<1\\
D_1(2)&=2\cos\frac{\pi}p\cos\frac{\pi}q\\
D_1(3)&=\cos\frac{\pi}p~[4(\cos\frac{\pi}q)^2-1]
\end{align*}
On dispose du résultat suivant:

\begin{prop}\label{p_2}
Soit $\Gamma(2,p,q)$ avec $p,q$ deux entiers tels que $p\leq q$ et
$\frac1p+\frac1q < \frac12$.
\begin{enumerate}
\item Si $p \geq 4$ et $q \geq 6$, les deux distances de
translation les plus petites parmi
les éléments hyperboliques de niveau $1$ sont, dans l'ordre croissant:\\
$l_1(2)=2\mathrm{Argcosh}~2~\cos\frac{\pi}p\cos\frac{\pi}q$ \\
$l_1(3)=2\mathrm{Argcosh}~\cos\frac{\pi}p~[4(\cos\frac{\pi}q)^2-1]$
\item Si $p=3$ où bien $(p,q)\in\{(4,5),(5,5) \}$, la plus petite
distance de translation parmi les éléments hyperboliques de niveau $1$ est:\\
$l_1(3)=2\mathrm{Argcosh}~\cos\frac{\pi}p~[4(\cos\frac{\pi}q)^2-1]$
\end{enumerate}
\end{prop}

\textbf{Preuve de la proposition \ref{p_2}:}\\
Un calcul rapide montre que $l_1(2)=l_1(3)$ si et seulement si
$q=5$. Il est de plus évident que $D_1(2)$ est strictement
supérieur à $1$ si $p\geq 4$ et strictement inférieur à $1$ dès
que $p=3$.$\square$

\subsection{Les hyperboliques de niveau 2}\label{s_32}

Soit $\gamma$ un élément hyperbolique de $\Gamma(2,p,q)$ déplaçant
au moins un sommet $s_0$ du pavage d'une distance combinatoire
$2$, c'est à dire $D(s_0,\gamma s_0)=2$. Cette isométrie se
décompose là encore en deux réflexions dont les axes sont des
géodésiques $G$ et $G'$ du pavage $\mathcal{P}_{0}$. Expliquons
comment s'en persuader.\\

Considérons $\beta=\{s_0x,x\gamma s_0\}$ un chemin de longueur $2$
reliant $s_0$ à son image. Posons maintenant $y=\gamma x$ et
considérons dès lors pour $G$ la bissectrice de l'angle
$\angle(xs_0,x\gamma s_0)$ et pour $G'$ la bissectrice de l'angle
$\angle(\gamma s_0x,\gamma s_0 y)$.

La longueur de translation de
$\gamma=\bar{r}_{G'_{k'}}\bar{r}_{G_k}$ est alors donnée par
\[l_2(k,k')=2\mathrm{Argcosh}~D_2(k,k')~;~ 1 \leq k,k' \leq q-1 \]
dès que la quantité
\[D_2(k,k')=\mid \sin\frac{k'\pi}q\sin\frac{k\pi}q\cosh 2c -\cos\frac{k\pi}q\cos\frac{k'\pi}q \mid~;~1\leq k,k'\leq q-1 \]
est strictement plus grande que $1$. Dans les formules
précédentes,
\[ \cosh2c=\frac{2(\cos\frac{\pi}p)^2+(\cos\frac{\pi}q)^2-1}{(\sin\frac{\pi}q)^2}\]

Il convient de remarquer que $D_2(q-1,k)$ est égal à $D_2(1,q-k)$
et cela pour tout $k$ compris entre $1$ et $q-1$.

Rangons ces valeurs dans l'ordre croissant. On dispose pour cela
du résultat élémentaire suivant:

\begin{lemme}\label{l_4}
Soit $K$ et $K'$ deux entiers compris entre $1$ et la partie
entière de $q/2$. Alors pour tout $k$ compris entre $K$ et $q-K$
et tout $k'$ compris entre $K'$ et $q-K'$, on a \[
\sin\frac{k\pi}q\sin\frac{k'\pi}q\cosh2c-\cos\frac{k\pi}q\cos\frac{k'\pi}q
\geq
\sin\frac{K\pi}q\sin\frac{K'\pi}q\cosh2c-\cos\frac{K\pi}q\cos\frac{K'\pi}q\]
\end{lemme}

\textbf{Preuve du lemme \ref{l_4}:}\\
On constate que pour les valeurs de $k$ indiquées,\[
\sin\frac{k\pi}q\geq\sin\frac{K\pi}q>0 \] et qu'il en est de même
pour $k'$ et $K'$. On a donc
\[
\sin\frac{k\pi}q\sin\frac{k'\pi}q\geq\sin\frac{K\pi}q\sin\frac{K'\pi}q>0
\] Maintenant,
\[ \cos\frac{k\pi}q\leq\cos\frac{K\pi}q \] avec le même résultat pour $k'$
et $K'$. On en déduit
\[ \cos\frac{k\pi}q \cos\frac{k'\pi}q\leq\cos\frac{K\pi}q\cos\frac{K'\pi}q \]
Ceci est évident sans se préoccuper du signe du terme de gauche
car le terme de droite est toujours positif par hypothèse sur $K$
et $K'$. On conclut en remarquant que $\cosh 2c$ est toujours
strictement positif. $\square$\\

Nous énonçons maintenant le principal résultat de ce paragraphe:
\begin{prop}\label{p_3}
Soit $\Gamma(2,p,q)$ avec $p,q$ deux entiers tels que $p\leq q$ et
$\frac1p+\frac1q < \frac12$.
\begin{enumerate}
\item Si $p \geq 5$, les deux distances de translation les plus
petites parmi
les éléments hyperboliques de niveau $2$ sont, dans l'ordre croissant:\\
$l_2(1,2)=2\mathrm{Argcosh}~[\cos\frac{\pi}q(4(\cos\frac{\pi}p)^2-1)]$\\
$l_2(1,q-1)=2\mathrm{Argcosh}~[
2(\cos\frac{\pi}p)^2+2(\cos\frac{\pi}q)^2-1]$ \item Si $p=4$ la
plus petite distance de translation parmi
les éléments hyperboliques de niveau $2$ est:\\
$l_2(1,q-1)=2\mathrm{Argcosh}~[
2(\cos\frac{\pi}p)^2+2(\cos\frac{\pi}q)^2-1]$ \item Si $p=3$ et
$q=7$,la plus petite distance de translation parmi
les éléments hyperboliques de niveau $2$ est:\\
$l_2(1,q-1)=2\mathrm{Argcosh}~[2(\cos\frac{\pi}q)^2-\frac12]$
\item Si $p=3$ et $q\geq 8$, les deux plus petites distances de
translation parmi
les éléments hyperboliques de niveau $2$ sont:\\
$l_2(1,q-1)=2\mathrm{Argcosh}~[2(\cos\frac{\pi}q)^2-\frac12]$\\
$l_2(1,6)=2\mathrm{Argcosh}~[2\cos\frac{\pi}q(2(\cos\frac{\pi}q)^2-1)]$
\end{enumerate}
\end{prop}

\textbf{Preuve de la proposition \ref{p_3}:}\\
\underline{\textit{Un calcul préliminaire}}\\
Pour clarifier la démonstration, nous commençons par obtenir une
première minoration des quantités $D_2(k,k')$ en distinguant deux
cas ($p\geq 4$ et $p=3$).

Si $p\geq 4$, la quantité entre les valeurs absolues est positive
dès que $1 \leq k \leq q-1$ et $2\leq k' \leq q-2$. Cela résulte
d'un calcul élémentaire et de l'utilisation  du lemme \ref{l_4}.
Ceci permet de faciliter la manipulation des quantités $D_2(k,k')$
en supprimant les valeurs absolues dans leur expression initiale.
On a alors, pour tout $2\leq k,k'\leq q-2$,
\begin{align*}
D_2(k,k')&=\sin\frac{k'\pi}q\sin\frac{k\pi}q\cosh 2c
-\cos\frac{k\pi}q\cos\frac{k'\pi}q \\
& \geq\sin\frac{2\pi}q\sin\frac{2\pi}q\cosh 2c
-\cos\frac{2\pi}q\cos\frac{2\pi}q\\
& =D_2(2,2)
\end{align*}
ceci résultant du lemme \ref{l_4}.

Si $p=3$, la quantité entre les valeurs absolues définissant
$D_2(k,k')$ est positive dès que $1 \leq k \leq q-1$ et $3\leq k'
\leq q-3$ et on a de même, pour tout $2\leq k \leq q-2$ et tout $3
\leq k' \leq q-3$,
\begin{align*}
D_2(k,k')&=\sin\frac{k'\pi}q\sin\frac{k\pi}q\cosh 2c
-\cos\frac{k\pi}q\cos\frac{k'\pi}q \\
& \geq\sin\frac{2\pi}q\sin\frac{3\pi}q\cosh 2c
-\cos\frac{2\pi}q\cos\frac{3\pi}q\\
& =D_2(2,3)
\end{align*}
Il convient donc d'étudier séparemment le cas $p=3$.\\

\underline{\textit{Le~cas~$p\geq 4$}}\\
Montrons maintenant la proposition dans le cas où $p$ est
supérieur ou égal à $4$.

On constate que
\[ D_2(1,2) <D_2(1,q-1)\leq D_2(1,3)\leq D_2(2,2) \]
avec un cas d'égalité qui caractérise la situation $p=4$.

En effet, on dispose des formules
\begin{align*}
&D_2(1,2)=\cos\frac{\pi}q[4(\cos\frac{\pi}p)^2-1]\\
&D_2(1,q-1)=2(\cos\frac{\pi}q)^2+2(\cos\frac{\pi}p)^2-1\\
&D_2(1,3)=8(\cos\frac{\pi}q)^2(\cos\frac{\pi}p)^2-2(\cos\frac{\pi}p)^2-2(\cos\frac{\pi}q)^2+1\\
&D_2(2,2)=8(\cos\frac{\pi}q)^2(\cos\frac{\pi}p)^2-1
\end{align*}
et le résultat provient des constations suivantes
\begin{align*}
&D_2(1,3)-D_2(1,q-1)=2(2(\cos\frac{\pi}q)^2-1)(2(\cos\frac{\pi}p)^2-1)\geq 0\\
&D_2(1,q-1)-D_2(1,2)> 2(\cos\frac{\pi}p)^2-(4(\cos\frac{\pi}p)^2-1)+2(\cos\frac{\pi}p)^2-1=0\\
&D_2(2,2)-D_2(1,3)\geq 4(\cos\frac{\pi}p)^2-2\geq0
\end{align*}
La première différence est nulle si et seulement si $p=4$ car on
demande $p\leq q$. Remarquons ensuite que si $p$ est supérieur ou
égal à $5$, $D_2(1,2)>1$ alors que $D_2(1,2)<1<D_2(1,q-1)$ si $p$
prend la valeur $4$.

Il reste à démontrer qu'aucune valeur $D_2(k,k')$ ne vient
s'intercaler dans la liste proposée. Si $k$ et $k'$ sont entre $2$
et $q-2$, cela résulte de la minoration obtenue en début de
preuve. Si $k=1$, on constate que $D_2(1,k)$ est supérieur à
$D_2(1,3)$ dès que $k$ est compris entre $3$ et $q-3$, toujours en
utilisant le lemme \ref{l_4}. Clairement, $D_2(1,1)<1$. Un calcul
montre ensuite que, si $p\geq 4$,
\begin{align*}
D_2(1,q-2)-D_2(1,q-1)&=4(\cos\frac{\pi}p)^2\cos\frac{\pi}q+4(\cos\frac{\pi}q)^3-3\cos\frac{\pi}q\\
&-2(\cos\frac{\pi}q)^2-2(\cos\frac{\pi}p)^2+1\\
&=
2(2\cos\frac{\pi}q-1)((\cos\frac{\pi}p)^2+(\cos\frac{\pi}q)^2)+1-3\cos\frac{\pi}q\\&\geq0
\end{align*}
la dernière égalité résultant d'une majoration directe si $p\geq4$
et $q\geq5$. Enfin, si $k=q-1$, on utilise que
$D_2(q-1,k')=D_2(1,q-k')$. Ceci achève la situation où
$p$ est supérieur à $4$.\\

\underline{\textit{Le cas $p=3$}}\\
Considérons maintenant le cas où $p=3$. On montre alors que
\[1<D_2(1,q-1) \leq D_2(1,6) \] avec un cas d'égalité qui
caractérise  la situation où $q=7$. Pour cela, on examine les
formules
\begin{align*}
D_2&(1,q-1)=D_2(1,5)=2(\cos\frac{\pi}q)^2-\frac12\\
D_2&(1,6)=D_2(2,q-1)=D_2(2,3)=2(\cos\frac{\pi}q)(2(\cos\frac{\pi}q)^2-1)\\
\end{align*}
On constate ensuite que $D_2(1,k)$ est supérieur à $D_2(1,6)$ dès
que $k$ est compris entre $6$ et $q-6$ puis on examine les valeurs
restantes comme précédemment pour confirmer qu'elles ne viennent
par s'intercaler dans la suite précédente. $\square$

\section{Aspect pratique de l'étude: deuxième partie} \label{s_4}

Pour $p=3$, on démontre (section \ref{s_422}) que tout élément
hyperbolique de $\Gamma(2,p,q)$ ayant une distance de translation
inférieure ou égale à $l_1(4)$ est de niveau $1$ ou $2$. Autrement
dit, en appliquant la propriété \ref{p_1} si $l_0=l_1(4)$, on
obtient $n_0=2$. Rappelons que l'on a étudié ces éléments dans la
section précédente.

Pour $p\in\{4,5\}$ et $p\geq 11$, la section \ref{s_421} montre
qu'il en est de même pour les éléments hyperboliques de
$\Gamma(2,p,q)$ ayant une distance de translation inférieure ou
égale à $l_2(1,q-1)$.

Enfin, pour $p\in\{6,7,8,9,10\}$, la section \ref{s_423} explique
pourquoi tout hyperbolique de $\Gamma(2,p,q)$ ayant une distance
de translation inférieure ou égale à $l_1(2)$ est également de
niveau $1$ ou $2$.

Dans chaque cas, l'objectif est donc d'expliquer pourquoi un
élément hyperbolique ayant une telle distance de translation ne
peut pas déplacer un élément du pavage d'une distance combinatoire
supérieure ou égale à $3$ et il est donc nécessaire, en
préliminaire, d'examiner de plus près la quantité $\rho(3)$.

\subsection{Calcul de $\rho(3)$} \label{s_41}

Fixons $s_0$ un sommet du pavage et considérons l'ensemble des
sommets qui constituent la sphère de centre $s_0$ et de rayon $3$
pour la distance combinatoire. Chaque point de cette sphère sera
écrit sous la forme $s(k_1,k_2)$ avec $k_1$ et $k_2$ entre $1$ et
$q-1$. Cherchons celui qui est géométriquement (i.e. au sens de la
distance hyperbolique $d$) le plus proche de $s_0$.\\

Començons par établir le

\begin{lemme}\label{l_5}
Soit $s$ et $s_0$ deux sommets de $\mathcal{P}$ tels que
$D(s,s_0)=3$. Si $\beta=\{s_0x,xy,ys\}$ est un chemin minimal de
$\mathcal{P}$ reliant $s$ à $s_0$, alors la mesure géométrique de
l'angle $\angle(yx,ys_0)$ est inférieure ou égale à $\pi/q$.
\end{lemme}

En reprenant la mesure de l'angle orienté, le lemme affirme que
\[\angle(yx,ys_0)\in[0,\frac{\pi}q]\cup[2\pi-\frac{\pi}q,\frac{\pi}q] \]

\textbf{Preuve du lemme \ref{l_5}:}\\
Commençons par signaler le fait que dans le triangle isocèle
$s_0xy$, l'angle à la base $\angle(yx,yx_0)$ est une fonction
décroissante de l'angle au sommet $\angle(xs_0,xy)$. Cela résulte
des formules classiques de
trigonométrie hyperbolique.\\

\underline{\textit{Le cas $p\geq 4$}}\\
Ici, on écrit $s=s(k_1,k_2)$ avec $k_1$ et $k_2$ entre $1$ et
$q-1$. La mesure géométrique de l'angle $\angle(xs_0,xy)$ est donc
minimale pour $k_1=1$, et un simple calcul (on pourra dessiner un
$p$-gone et représenter les points $s_0,x,y$) montre que si $p\geq
4$, l'angle $\angle(yx,ys_0)$ a dans ce cas une mesure géométrique
toujours
inférieure ou égale à $\pi/q$ avec comme cas d'égalité la valeur $p=4$.\\

\underline{\textit{Le cas $p=3$}}\\
On écrit $s=s(k_1,k_2)$ avec $k_1$ et $k_2$ entre $2$ et $q-2$.
Remarquons que si $k_1=1$ ou $k_1=q-1$, l'angle $\angle(yx,ys_0)$
a une mesure géométrique égale à $2\pi/q$ et le raisonnement
précédent n'est a priori pas exploitable. Cependant, avec
l'hypothèse où $k_2$ doit être entre $2$ et $q-2$, la mesure
géométrique de l'angle $\angle(xs_0,xy)$ est minimale pour $k_1=2$
ou $k_1=q-2$. On trouve alors une mesure géométrique égale à
$\pi/q$ pour l'angle $\angle(xs_0,xy)$ (on pourra là encore
dessiner un $p$-gone pour s'en convaincre). Comme il s'agit du cas
extrémal,
cela conclut la preuve. $\square$\\

On déduit de cela le résultat suivant

\begin{lemme}\label{l_6}
Supposons $p\geq 3$. Soit $K$ un entier compris entre $1$ et la
partie entière de $q/2$. Pour tout $k_1$ entre $1$ et $q-1$ et
$k_2$ entre $K$ et $q-K$,
\[ d(k_1,k_2)\geq \mathrm{Min}\{d(k_1,K),d(k_1,q-K)\}
\]
\end{lemme}

\textbf{Preuve du lemme \ref{l_6}:}:\\
On se fixe $s_0$ un point du pavage et un couple d'entier
$(k_1,k_2)$ auquel on associe un point $s=s(k_1,k_2)$ de la sphère
combinatoire de rayon $3$ et de centre $s_0$. On dispose d'un
chemin $\beta=\{s_0x,xy,ys\}$ associé à ce couple d'entier. Dans
le triangle $s_0ys$, la longueur $s_0s$ est minimale quand la
mesure géométrique de l'angle $\angle(ys_0,ys)$ est minimal. En
notant $\alpha$ la mesure géométrique de l'angle
$\angle(yx,ys_0)$, on peut distinguer deux cas:

Si $2k_2\pi/q$ est dans $]0,\pi]$, la mesure géométrique de
$\angle(ys_0,ys)$ est égale à \[2k_2\pi/q-\alpha\] Elle est donc
minimale pour $k_2=K$.

Si $2k_2\pi/q$ est dans $]\pi,2\pi]$, cette mesure géométrique est
égale à \[2(q-k_2)\pi/q + \alpha\] et est minimale pour $k_2=q-K$.

Ceci
termine la preuve.$\square$\\

Ce lemme affirme par symétrie que
\[ d(k_1,k_2)\geq
\mathrm{Min}\{d(K,k_2),d(q-K,k_2)\} \] si $k_1$ varie entre $K$ et
$q-K$ et $k_2$ est un entier quelconque entre $1$ et $q-1$.

On déduit du lemme \ref{l_6} la propriété suivante.

\begin{prop}\label{p_5}La valeur de $\rho(3)$ est déterminée de la
manière suivante:
\begin{enumerate}
\item si $p\geq 6$, $\rho(3)=d(s_0,s(1,1))$\item si $p=4$,
$\rho(3)=d(s_0,s(1,2))$ \item si $p=5$, $\rho(3)=d(s_0,s(1,q-1))$
\item si $p=3$, $\rho(3)=d(s_0,s(2,2))$
\end{enumerate}
\end{prop}

\textbf{Preuve de la proposition \ref{p_5}:}\\
On se fixe $s_0$ un sommet du pavage et on désigne les points de
la sphère de centre $s_0$ et de rayon $3$ pour la distance
combinatoire par les sommets $s(k_1,k_2)$ avec $k_1$ et $k_2$
compris entre $1$ et $q-1$.\\

\underline{\textit{Le cas $p\geq6$}}\\
Dans ce cas, pour tous $k_1$ et $k_2$ compris entre $1$ et $q-1$,
le point $s(k_1,k_2)$ est à une distance combinatoire $3$ de
$s_0$.

Pour tout $k_1$ fixé, en appliquant le lemme \ref{l_6} à deux
reprises, on a
\[ d(k_1,k_2)\geq
\mathrm{Min}\{d(1,1),d(1,q-1)\} \]

Montrons que, pour $p\geq 6$, $d(1,q-1)$ est toujours supérieur à
$d(1,1)$.

En écrivant $\beta_1=\{s_0x,xy,ys(1,1)\}$ et
$\beta_2=\{s_0x,xy,ys(1,q-1)\}$ des chemins du pavage associés à
$s(1,1)$ et $s(1,q-1)$, on raisonne dans les triangles
$s_0ys(1,1)$ et $s_0ys(1,q-1)$. Notons $\alpha$ la mesure
géométrique de l'angle $\angle(yx,ys_0)$. L'angle en $y$ du
premier triangle (qui est $2\pi/q-\alpha$) est strictement
inférieur à celui du second (qui est $2\pi/q+\alpha$). On termine
en invoquant le fait que dans un triangle hyperbolique, la
longueur d'un
côté est une fonction croissante de l'angle opposé à ce côté.\\

\underline{\textit{Le cas $p=4$ ou $p=5$}}\\
Ici, si $k_1=1$ et $k_2=1$, on a $D(s_0,s(1,1))<3$ donc il s'agit
de retirer cette possibilité.

Le raisonnement précédent, s'appuyant sur le lemme \ref{l_6},
permet d'établir, d'une part, que pour tous $k_1,k_2$ entre $2$ et
$q-2$,
\[d(k_1,k_))\geq
\mathrm{Min}\{d(2,2),d(2,q-2)\} \] et d'autre part que pour $k_1$
entre $1$ et $q-1$ et $k_2$ entre $2$ et $q-2$,
\[d(k_1,k_2)\geq
\mathrm{Min}\{d(1,2),d(q-1,2)\} \] Une nouvelle application permet
d'affirmer
\[d(q-1,2)\geq d(q-1,1) \] en remarquant que
$s(q-1,1)$ est bien à une distance combinatoire $3$ de $s_0$ ici.
Reste donc à comparer $d(1,2)$ et $d(1,q-1)$ pour $p=4$ et $p=5$
pour conclure.

En écrivant $\beta_1=\{s_0x,xy,ys(1,2)\}$ et
$\beta_2=\{s_0x,xy,ys(1,q-1)\}$ des chemins du pavage associés à
$s(1,2)$ et $s(1,q-1)$, on raisonne dans les triangles
$s_0ys(1,2)$ et $s_0ys(1,q-1)$. Les angles au sommet $y$ sont
égaux dans ces deux triangles si $p=4$ alors que l'angle en $y$ du
triangle $s_0ys(1,2)$ est plus grand que celui du triangle
$s_0ys(1,q-1)$ si $p=5$ (on pourra dessiner un $p$-gone pour s'en
convaincre). On conclut comme précédemment en invoquant le fait
que, dans un triangle hyperbolique, la longueur d'un
côté est une fonction croissante de l'angle opposé à ce côté.\\

\underline{\textit{Le cas $p=3$}}\\
Dans ce dernier cas, $k_1$ et $k_2$ doivent évoluer entre $2$ et
$q-2$ pour que l'on obtienne bien un point $s(k_1,k_2)$ qui soit à
une distance $3$ de $s_0$. On utilise  le lemme \ref{l_6} pour
établir que
\[d(s_0,s(k_1,k_2))\geq
\mathrm{Min}\{d(2,2),d(2,q-2)\} \] et il reste à comparer ces deux
longueurs pour conclure: il s'agit du même argument que lorsqu'on
voulait comparer $d(1,1)$ et $d(1,q-1)$ dans le cas où $p\geq 6$.

Ceci termine la preuve du lemme. $\square$\\

Reste à calculer explicitement ces valeurs en fonction de $p$ et
$q$. Pour cela, on utilise le changement de variable suivant:
\[ X=\cos\frac{\pi}p~;~Y=\cos\frac{\pi}q \]
On démontre alors la proposition

\begin{prop}\label{p_6}
Avec les notations précédentes, on dispose des formules suivantes:
\begin{align*}
\cosh d(1,2)&=\frac{32X^2Y^2(X^2+Y^2-1)(4X^2-1)+2X^2+Y^2-1}{1-Y^2}\\
\cosh d(1,1)&=\frac{16X^2(X^2+Y^2-1)(2X^2-1)+2X^2+Y^2-1}{1-Y^2}\\
\cosh d(1,q-1)&=\frac{(2X^2+Y^2-1)[16X^2(X^2+Y^2-1)+1]}{1-Y^2}\\
\cosh d(2,2)&=\frac{64X^2Y^2(X^2+Y^2-1)(8Y^2X^2-1)+2X^2+Y^2-1}{1-Y^2}\\
\end{align*}
\end{prop}

\textbf{Preuve de la proposition \ref{p_6}:}\\
Dans toute cette preuve, $s_0$ est un sommet fixé du pavage
$\mathcal{P}$.\\

\underline{\textit{Une première formule générale}}\\
On suppose que, $2k_1\pi/q$ et $2k_2\pi/q$ sont strictement
inférieurs à $\pi$. On montre ici que
\begin{align*} \cosh
d(k_1,k_2)&=\cosh2c(\sinh2c)^2(1-\cos\frac{2k_1\pi}q)(1-\cos\frac{2k_2\pi}q)\\
&-(\sinh 2c)^2\sin\frac{2k_1\pi}q\sin\frac{2k_2\pi}q+\cosh 2c
\end{align*}

Si $s=s(k_1,k_2)$, notons $s_0,x,y,s$ les sommets consécutifs du
quadrilatère hyperbolique formé par le chemin reliant $s_0$ à $s$
dans $\mathcal{P}$ en y rajoutant le segment géodésique de
$\mathbf{H}$ reliant $s$ à $s_0$. On a
\begin{align*}
\cosh d(s_0,s) &=\cosh 2c \cosh d(s_0,y)-\sinh 2c\sinh
d(s_0,y)\cos(\frac{2k_2\pi}q-\angle(yx,ys_0))\\
&=\cosh 2c \cosh d(s_0,y)-\sinh 2c\sinh
d(s_0,y)\cos\frac{2k_2\pi}q\cos\angle(yx,ys_0)\\
&-\sinh 2c\sinh d(s_0,y)\sin\frac{2k_2\pi}q\sin\angle(yx,ys_0)
\end{align*}
or
\begin{align*}
\cosh d(s_0,y)&=(\sinh 2c)^2(1-\cos\frac{2k_2\pi}q)+1~;\\
\sin\angle(yx,ys_0)\sinh d(s_0,y)&=\sinh 2c \sin\frac{2k_1\pi}q
\end{align*}
donc
\begin{align*}
\cosh d(s_0,s)&=\cosh 2c [(\sinh
2c)^2(1-\cos\frac{2k_1\pi}q)+1]\\
&-\sinh 2c \sinh d(s_0,y)
\cos\frac{2k_2\pi}q\cos\angle(yx,ys_0)\\
&-(\sinh 2c)^2\sin\frac{2k_1\pi}q\sin\frac{2k_2\pi}q
\end{align*}
mais on sait que

\[ \cosh 2c =\cosh 2c\cosh d(s_0,y)-\sinh 2c\sinh
d(s_0,y)\cos\angle(yx,ys_0)\]ce qui entraîne \[\sinh 2c\sinh
d(s_0,y)\cos\angle(yx,ys_0)=\cosh 2c [\cosh d(s_0,y)-1]
\]
d'où
\begin{align*}
\cosh d(s_0,s)&=\cosh 2c [(\sinh
2c)^2(1-\cos\frac{2k_1\pi}q)+1]\\
&-(\sinh 2c)^2 \cosh 2c\cos\frac{2k_2\pi}q(1-\cos\frac{2k_1\pi}q)\\
&-(\sinh 2c)^2\sin\frac{2k_1\pi}q\sin\frac{2k_2\pi}q\\
&=\cosh 2c(\sinh 2c)^2(1-\cos\frac{2k_1\pi}q)+\cosh 2c\\
&-(\sinh 2c)^2\sin\frac{2k_1\pi}q\sin\frac{2k_2\pi}q\\
&-\cosh 2c(\sinh 2c)^2(1-\cos\frac{2k_1\pi}q)\cos\frac{2k_2\pi}q\\
&=\cosh2c(\sinh2c)^2(1-\cos\frac{2k_1\pi}q)(1-\cos\frac{2k_2\pi}q)\\
&-(\sinh 2c)^2\sin\frac{2k_1\pi}q\sin\frac{2k_2\pi}q+\cosh 2c
\end{align*}
ce qui termine la preuve .\\

\underline{\textit{Applications de la formule obtenue}}\\
La formule ci-dessus permet de trouver toutes les quantités
cherchées à l'exception de $d(1,q-1)$. Examinons en détail le cas
de $d(1,2)$. Via les formules
\begin{align*}
& \cosh 2c =2(\cosh c)^2-1=\frac{2X^2+Y^2-1}{1-Y^2}~;\\
& (\cosh 2c)^2-1=4\frac{X^2(X^2+Y^2-1)}{(1-Y^2)^2}~;\\
& \cos \frac{2\pi}q=2Y^2-1~;\\
& \cos\frac{4\pi}q= -8Y^2(1-Y^2)+1~;\\
& \sin \frac {2\pi}q\sin\frac{4\pi}q=8Y^2(1-Y^2)(2Y^2-1)
\end{align*}
qu'il convient de faire intervenir dans la formule générale
obtenue plus haut, on obtient
\begin{align*}
\cosh d(1,2) &= \cosh 2c((\cosh 2c)^2-1)(1-\cos\frac{2\pi}q)(1-\cos\frac{4\pi}q)\\
&-((\cosh 2c)^2-1)\sin\frac{2\pi}q\sin\frac{4\pi}q + \cosh 2c\\
&=\frac{(2X^2+Y^2-1)4X^2(X^2+Y^2-1)(2-2Y^2)8Y^2(1-Y^2)}{(1-Y^2)^3}\\
&-\frac{4X^2(X^2+Y^2-1)8Y^2(1-Y^2)(2Y^2-1)}{(1-Y^2)^2}+\frac{2X^2+Y^2-1}{1-Y^2}\\
&=\frac{32(X^2+Y^2-1)X^2Y^2(4X^2-1)+2X^2+Y^2-1}{1-Y^2}
\end{align*}
C'est ce que l'on voulait. On laisse au lecteur le soin de
détailler les autres cas.\\

\underline{\textit{Calcul de $d(1,q-1)$}}\\
De la même manière que ci-dessus, on note $s_0,x,y,s$ les sommets
consécutifs du pavage associés au chemin de longueur combinatoire
$3$ reliant $s_0$ à $s=s(1,q-1)$. On a
\begin{align*}
\cosh d(s_0,s) &=\cosh 2c \cosh d(s_0,y)-\sinh 2c\sinh
d(s_0,y)\cos(\frac{2\pi}q+\angle(yx,ys_0))\\
&=\cosh 2c \cosh d(s_0,y)-\sinh 2c\sinh
d(s_0,y)\cos\frac{2\pi}q\cos\angle(yx,ys_0)\\
&+\sinh 2c\sinh d(s_0,y)\sin\frac{2\pi}q\sin\angle(yx,ys_0)
\end{align*}
le calcul est essentiellement identique à celui fait dans le cas
où les angles étaient inférieurs à $\pi$.
On trouve alors le résulat annoncé. $\square$\\

Avant, de poursuivre, signalons qu'il existe une méthode plus
rapide pour calculer $d(1,1)$.

En désignant par $O$ le centre du polygone du pavage $\mathcal{P}$
admettant $s_0$ et $s$ parmi ses sommets, et en appelant $b$ la
distance hyperbolique de $O$ à une arète quelconque de ce
polygone, on a
\[ \cosh b = \frac{\cos\frac{\pi}q}{\sin{\frac{\pi}p}} \]
et
\[ \cosh d(s_0,s(1,1))=(\sinh b)^2(1-\cos\frac{6\pi}p)+1 \]

\subsection{Recherche des longueurs les plus courtes} \label{s_42}
On va maintenant utiliser le calcul de $\rho(3)$ pour appliquer la
proposition \ref{p_1} et décrire les distances de translation du
spectre des longueurs de $\Gamma(2,p,q)$ qui se trouvent être
inférieures à une certaine valeur $l_0$. Cette longueur $l_0$
étant différente suivant les valeurs de $p$, nous exposons les
résultats de cette partie en fonction du paramètre $p$.

\subsubsection{Le cas où $p=4,5$ ou $p\geq 11$} \label{s_421}

Supposons dans ce paragraphe que $p\in\{4,5\}$ ou bien $p\geq 11$.

On montre que les éléments hyperboliques $\gamma$ de
$\Gamma(2,p,q)$ qui sont de longueur inférieure ou égale à
$l_0=l_2(1,q-1)$ on un niveau $\lambda(\gamma)$ nécessairement
inférieur ou égal à $2$.

Nous allons pour cela mettre en pratique la proposition \ref{p_1}:
il suffit de prouver que
\[ \rho(3)>C(l_2(1,q-1))\]

C'est l'objet de la

\begin{prop}\label{p_7}
Si $p\in\{4,5\}$ ou $p\geq 11$, tout élément hyperbolique $\gamma$
de $\Gamma(2,p,q)$ tel que $\lambda(\gamma)\geq 3$ a une distance
de translation strictement supérieure à $l_2(1,q-1)$.
\end{prop}

\textbf{Preuve de la proposition \ref{p_7}:}\\
Examinons ce qui se passe pour les différentes valeurs de $p$.\\

\underline{\textit{Le cas $p=4$}}\\
Si $p=4$ on sait par la proposition \ref{p_5} que $\rho(3)=d(2,1)$
et il suffit d'établir que
\[ \cosh d(2,1)>(\cosh c)^2(\cosh l_2(1,q-1) -1)+1 \]

Le terme de gauche a été calculé à la proposition \ref{p_6}. Le
terme de droite est donné par
\begin{align*}
(\cosh c)^2(\cosh l_2(1,q-1)-1)+1&=2((D_2(1,q-1))^2-1)\frac{X^2}{1-Y^2}+1\\
&=\frac{2X^2((2X^2+2Y^2-1)^2-1)+1-Y^2}{1-Y^2}\\
&=\frac{2X^2(2X^2+2Y^2-2)(2X^2+2Y^2-1+1)+1-Y^2}{1-Y^2}\\
&=\frac{8X^2(X^2+Y^2-1)(X^2+Y^2)+1-Y^2}{1-Y^2} \end{align*}

Ce qui mène à l'étude de la fonction
\begin{align*}
h(X,Y)
&=\frac{2(X^2+Y^2-1)~[4X^2(16X^2Y^2-5Y^2-X^2)+1]}{1-Y^2}\\
&=\frac{g(X,Y)}{1-Y^2}
\end{align*}
et l'on conclut en remarquant que $g(\sqrt{2}/2,Y)$ est
strictement positive sur $[\cos\pi/5,1[$.

En effet, on observe que
\[4X^2(16X^2Y^2-5Y^2-X^2)+1=6Y^2>0 \]
ce qui termine la preuve dans le cas $p=4$.\\

\underline{\textit{Le cas $p=5$}}\\
Cette fois ci, $\rho(3)=d(1,q-1)$ et il s'agit de montrer que
\[ \cosh d(1,q-1)>(\cosh c)^2(\cosh l_2(1,q-1)-1)+1 \]
On constate pour cela que la fonction
\begin{multline*}
g(X,Y)=(2X^2+Y^2-1)(16X^2(X^2+Y^2-1)+1)
\\-8X^2(X^2+Y^2-1)(X^2+Y^2)-1+Y^2
\end{multline*}
est strictement positive pour $X=\cos\pi/5$ sur
$[\cos\frac{\pi}5,1[$ ce qui termine le cas $p=5$.\\

\underline{\textit{Le cas $p\geq 11$}}\\
On a $\rho(3)=d(1,1)$. Il suffit de montrer
\[ \cosh d(1,1) >(\cosh
c)^2(\cosh l_2(1,q-1) -1) +1\]

qui est une conséquence du fait que la fonction
\begin{align*}
g(X,Y)&=(X^2+Y^2-1)~[~16~X^2(2X^2-1)-8~X^2(X^2+Y^2)+2~]
\\&=(X^2+Y^2-1)~h(X,Y)
\end{align*}
soit strictement positive  sur $[ \cos\frac{\pi}{11},1 [^2$ en
remarquant que sur ce produit,
\[h(X,Y)
> h(X,1)=24~X^4-24~X^2+2 > 0 \]
Ceci termine la preuve de la proposition \ref{p_7}.$\square$\\

On a donc déterminé toutes les valeurs du spectre des longueurs
des groupes $\Gamma(2,4,q),\Gamma(2,5,q)$ et
 $\Gamma(2,p,q),p\geq 11$ qui sont inférieures à $l_2(1,q-1)$.

\subsubsection{Le cas $p=3$} \label{s_422}

On étudie ici les groupes $\Gamma(2,3,q)$ avec $q\geq 7$. Là
encore, il s'agit de mettre en pratique la propriété \ref{p_1},
mais cette fois-ci avec $l_0=l_1(4)$.

On montre la proposition suivante

\begin{prop}\label{p_8}
Si $p=3$, alors tout élément hyperbolique $\gamma$ de
$\Gamma(2,3,q)$ tel que $\lambda(\gamma)\geq 3$ a une distance de
translation strictement supérieure à $l_1(4)$.
\end{prop}

\textbf{Preuve de la proposition \ref{p_8}:}\\
Soit $s$ un sommet du pavage $\mathcal{P}$. Il suffit de démontrer
que
\[\rho(3)>(\cosh c)^2(\cosh
l_1(4)-1)+1\] On rappelle pour cela que $\rho(3)=d(2,2)$, quantité
dont on dispose déjà d'une formule à la proposition \ref{p_6}. Un
calcul établit en outre
\[
(\cosh c)^2(\cosh l_1(4)-1)+1=\frac{Y^2[2(2Y^2-1)^2+\frac12]}{1-Y^2}
\]
On conclut en considérant la différence et en constatant
\[
g(Y)=24Y^6-32Y^4+12Y^2-1
\]est strictement positive sur $[\cos\frac{\pi}7,1[$, ce qu'une
étude élémentaire confirme. $\square$\\

On a donc déterminé toutes les distances de translation
inférieures à $l_1(4)$ se trouvant dans le spectre des longueurs
des groupes $\Gamma(2,3,q)$.

\subsubsection{Le cas $6\leq p \leq 10$} \label{s_423}

On suppose $6\leq p\leq 10$ et on montre la propriété suivante, en
appliquant encore une fois la proposition \ref{p_1} dans ce cas
précis avec $l_0=l_1(2)$.

\begin{prop}\label{p_9}
Soit $6\leq p\leq 10$. Alors tout élément hyperbolique $\gamma$ de
$\Gamma(2,p,q)$ vérifiant $\lambda(\gamma)\geq 3$ a une distance
de translation strictement supérieure à $l_1(2)$.
\end{prop}

\textbf{Preuve de la proposition~ \ref{p_9}:}\\
On montre que pour les valeurs de $p$ indiquées,
\[\rho(3)>(\cosh c)^2(\cosh l_1(2)-1)+1\]
Rappelons que  $\rho(3)=d(1,1)$ a déjà été calculé à la propriété
\ref{p_6}. On établit de plus que
\[
(\cosh c)^2(\cosh l_{1}(2)-1)+1=\frac{2X^2
(4X^2Y^2-1)+(1-Y^2)}{1-Y^2}
\]
Et nous sommes amené à étudier le signe de la fonction
\begin{multline*}
g(X,Y)=\\
16X^2(2X^2-1)(X^2 +Y^2-1)+2X^2+2Y^2-2-2X^2(4X^2Y^2-1)
\end{multline*} qui s'avère être strictement positif pour
$X\in[\sqrt{3}/2,1[$ et $Y\in[X,1[$.

En effet,
\[ g(X,Y) \geq 16X^2(2X^2-1)^2+4X^2-2-2X^2(4X^2-1)>0 \]sur le
domaine de définition étudié.$\square$\\

La preuve montre que le résultat de la proposition \ref{p_9} est
en fait valable pour $p\geq 6$.

On a donc décrit toutes les distances de translation inférieures à
$l_1(2)$ pour les groupes $\Gamma(2,p,q)$ qui manquaient à
l'appel.

\section{Rigidité spectrale}\label{s_5}

Dans cette partie, nous décrivons les premières valeurs du spectre
des longueurs des groupes $\Gamma(2,p,q)$ dans le théorème
\ref{t_1}. Nous en déduisons en particulier une formule explicite
de la systole dans le corollaire \ref{c_1} avant d'utiliser ces
formules pour démontrer la rigidité spectrale de ces objets. On
commence par établir cette rigidité pour les groupes vérifiant
$p\geq 11$ avant de généraliser.

\subsection{Description du début du
spectre}\label{s_51} Nous pouvons désormais énoncer le théorème
central de l'étude.
\begin{theoreme}\label{t_1}
Soit $\Gamma(2,p,q)$ avec $p$ et $q$ deux entiers vérifiant $p\leq
q$ et $\frac1p+\frac1q<\frac12$. Alors le début du spectre des
longueurs est donné par les formules suivantes:\\
$\mathrm{Lsp}~\Gamma(2,3,7)=\{l_2(1,q-1)=\dots \}$\\
$\mathrm{Lsp}~\Gamma(2,3,q)=\{l_2(1,q-1)=\dots<l_1(4)\dots \}$ pour tout $q\geq 8$\\
$\mathrm{Lsp}~\Gamma(2,4,5)=\{l_1(2)=\dots<l_2(1,q-1)\dots \}$ \\
$\mathrm{Lsp}~\Gamma(2,4,q)=\{l_1(2)=\dots<l_1(3)=\dots<l_2(1,q-1)\dots \}$ pour $q=6,7$\\
$\mathrm{Lsp}~\Gamma(2,4,q)=\{l_1(2)=\dots<l_2(1,q-1)\dots \}$ pour tout $q\geq8$\\
$\mathrm{Lsp}~\Gamma(2,5,5)=\{l_1(2)=\dots<l_2(1,q-1)\dots \}$ \\
$\mathrm{Lsp}~\Gamma(2,5,q)=\{l_1(2)=\dots<l_2(1,2)\dots \}$ pour tout $q\geq 6$\\
$\mathrm{Lsp}~\Gamma(2,p,q)=\{l_1(2)=\dots \}$ pour tout $p\in\{6,7,8,9,10\}$\\
$\mathrm{Lsp}~\Gamma(2,p,q)=\{l_1(2)=\dots<l_2(1,2)\dots \}$ pour tout $p\geq 11$\\
\end{theoreme}

avec
\begin{align*}
&l_2(1,q-1)=2~\textrm{Argch}~[2(\cos\frac{\pi}p)^2+2(\cos\frac{\pi}q)^2-1]\\
&l_2(1,2)=2~\textrm{Argch}~[\cos\frac{\pi}q[4(\cos\frac{\pi}p)^2-1]]\\
&l_1(2)=2~\textrm{Argch}~[2\cos\frac{\pi}p\cos\frac{\pi}q]\\
&l_1(3)=2~\textrm{Argch}~[\cos\frac{\pi}p[4(\cos\frac{\pi}q)^2-1]]
\end{align*}

\textbf{Preuve~du~théorème~\ref{t_1}}:\\
Nous devons ranger dans l'ordre croissant les valeurs trouvées
dans les propositions \ref{p_2} et \ref{p_3}.\\

\underline{\textit{Le cas $p=3$}}\\
En utilisant les formules données, on établit que
\[ 1<D_1(3)=D_2(1,q-1)\leq D_1(4)=D_2(1,6) \]
avec égalité si et seulement si $q=7$.\\

\underline{\textit{Le cas $p=4$:}}\\
Cette fois ci, on doit différencier deux cas.

Si $q=5,6,7$, on a
\[ D_2(1,2)<1<D_1(2)\leq D_1(3)< D_2(1,q-1) \]
avec égalité si et seulement si $q=5$.

Si par contre $q\geq8$, on constate que
\[ D_2(1,2)<1<D_1(2)<D_2(1,q-1)<D_1(3) \]
ce qui achève la preuve dans ce cas.\\

\underline{\textit{Le cas $p\geq5$:}}\\
On constate ici que
\[ 1<D_1(2)\leq D_2(1,2)\leq D_1(3) \]
avec égalité si et seulement si $q=5$. Toujours si $p=q=5$, on constate enfin que
\[D_1(3)<D_2(1,q-1) \]
ce qui achève la preuve.$\square$\\

En particulier, on obtient le

\begin{corollaire}\label{c_1} La
systole de $\Gamma(2,p,q)$ est déterminée par l'alternative
suivante:
\begin{enumerate}
\item Si $p\geq 4$, il s'agit de
$2~\mathrm{Argcosh}~[2\cos\frac{\pi}p\cos\frac{\pi}q]$ \item Si
$p=3$, il s'agit de
$2~\mathrm{Argcosh}~[2(\cos\frac{\pi}q)^2-\frac12]$
\end{enumerate}
\end{corollaire}

\subsection{Interprétation géométrique} \label{s_43}
Sachant que les classes de conjugaison des éléments hyperboliques
de distance de translation $l_0$ de $\Gamma(2,p,q)$ donnent
naissance aux géodésiques fermées de longueur $l_0$ de la surface
à points coniques $\textbf{H}/\Gamma(2,p,q)$, il est naturel de
chercher à représenter les géodésiques obtenues ci-dessus, et
notamment la systole.

Autrement dit, on cherche quelle est la forme géométrique des
géodésiques les plus courtes sur les surfaces à points coniques
considérées.

Rappelons que ces quotients sont de genre $0$ et admettent trois
points coniques non équivalents qui correspondent aux trois
classes de conjugaison d'éléments elliptiques présents dans
$\Gamma(2,p,q)$.

\subsection{Rigidité dans le cas $p \geq 11$}\label{s_52}

Dans cette partie, on montre que si $\Gamma(2,p,q)$ et
$\Gamma(2,p',q')$ vérifient $p\geq 11$ et $p'\geq 11$ et ont le
même spectre des longueurs, alors ils sont isométriques.

\begin{prop}\label{c_2}
Si $\Gamma(2,p,q)$ et $\Gamma(2,p',q')$ sont deux groupes de
triangles ayant le même spectre des longueurs et vérifiant $p \geq
11$ et $p' \geq 11$, alors $p=p'$ et $q=q'$.
\end{prop}

\textbf{Preuve de la proposition \ref{c_2}}:\\
On rappelle que le début du spectre est alors donné par
\[\mathrm{Lsp}~\Gamma(2,p,q)=\{l_1(2)=\dots<l_2(1,2)\dots \}\]
et que l'on connaît donc les longueurs $l_1(2)$ et $l_2(1,2)$ si
le spectre est donné.

Notons $l_1,l_2$ les deux premières valeurs distinctes du spectre
et déterminons
\[ L_i=\cosh\frac{l_i}2~;~i=1,2 \]
Il suffit d'établir que le système
\begin{displaymath}
\left\{
\begin{array}{ll}
2XY&=L_1\\
Y(4X^2-1)&=L_2
\end{array}\right.
\end{displaymath}
admet un unique couple de solution $(X,Y)$ dans $[0,1]^2$. Cela
revient à montrer que le polynôme
\[ P(u)=4L_1u^2-2L_2u-L_1 \]
a une unique racine dans $[0,1]$: cette racine est alors $X$. Le
polynôme $P$ admettant toujours deux racines réelles de signes
opposés, il en existe une seule dans l'intervalle
considéré. Ceci achève la preuve, modulo le changement de variable $X=\cos(\pi/p)$ et $Y=\cos(\pi/q)$. $\square$\\

\subsection{Le cas général} \label{s_53}

Avant d'entamer la preuve du théorème central, commençons par
énoncer et démontrer le résultat suivant:

\begin{lemme}\label{l4}
Un groupe $\Gamma(2,p',q')$ avec $p'\leq 8$ ne peut pas avoir le
même spectre des longueurs qu'un groupe $\Gamma(2,p,q)$ avec
$p\geq 9$.\\
Un groupe $\Gamma(2,3,q)$ ne peut pas avoir le même spectre des
longueurs qu'un groupe $\Gamma(2,p',q')$ avec $p'>5$.
\end{lemme}

\textbf{Preuve du lemme \ref{l4}:}\\
Commençons par montrer la première assertion.

Traitons d'abord le cas $p'\in\{4,5,6,7,8\}$. Dans le cas
d'isospectralité, on aurait égalité des systoles et donc\[
\cos\frac{\pi}p\cos\frac{\pi}q=\cos\frac{\pi}{p'}\cos\frac{\pi}{p'}<\cos\frac{\pi}8~;~p\geq
9
\]
Ceci impose \[ (\cos\frac{\pi}p)^2<\cos\frac{\pi}8 \] donc $p\leq
11$. On a de plus $q\leq 17$ en résolvant
\[\cos\frac{\pi}8>\cos\frac{\pi}p\cos\frac{\pi}q\geq\cos\frac{\pi}{9}\cos\frac{\pi}q\]
Il s'agit donc de trouver $4\leq p'\leq 8~;~q'\geq p'$ tels que
\[
\cos\frac{\pi}{p'}\cos\frac{\pi}{q'}=\cos\frac{\pi}{p}\cos\frac{\pi}{q}~;~9\leq
p\leq 11~;~p\leq q\leq 17
\]
En examinant séparemment les cas où $p'\in\{4,5,6,7,8\}$ on
constate qu'un tel $q'$ n'existe pas.

Montrons maintenant que $\Gamma(2,p,q)$ ne peut être isospectral à
un groupe $\Gamma(2,3,q')$. Pour cela, il suffit de constater que
ces groupes ne peuvent pas avoir la même systole car on aurait
\[\frac32>2(\cos\frac{\pi}{q'})^2-1=2\cos\frac{\pi}{p}\cos\frac{\pi}{q}
\]
ce qui montre que $p\leq 5$. La deuxième assertion est une
conséquence de ce qui vient d'être dit: la systole d'un groupe
$\Gamma(2,3,q)$ est nécessairement inférieure à $3/2$ ce qui est
impossible dès que $p$ est supérieur à $6$.
$\square$\\

On se propose d'établir maintenant le résultat principal de cet
article:

\begin{theoreme}\label{t_2}
Soit $\Gamma(2,p,q)$ et $\Gamma(2,p',q')$ deux groupes de
triangles ayant le même spectre des longueurs. Alors ces deux
groupes sont isométriques.
\end{theoreme}

\textbf{Preuve du théorème \ref{t_2}:}\\
Fixons nous un groupe $\Gamma(2,p,q)$ et montrons qu'il ne peut
pas être isospectral à un groupe $\Gamma(2,p',q')$ avec $p\neq p'$
ou bien $q\neq q'$.\\

\underline{\textit{Le cas $p=10$}}\\
On écarte les cas où $p'\leq 8$ en utilisant le lemme \ref{l4} .
Le cas $p'=9$ s'exclut par symétrie en considérant les calculs
faits pour $p=9$. Enfin, il n'est pas envisageable de trouver un
groupe $\Gamma(2,10,q)$ isospectral à un groupe $\Gamma(2,p',q')$
avec $p'\geq 11$ car cela imposerait, par égalité des systoles, \[
\cos\frac{\pi}{10}>\cos\frac{\pi}{10}\cos\frac{\pi}q=\cos\frac{\pi}{p'}\cos\frac{\pi}{q'}\geq
(\cos\frac{\pi}{p'})^2 \] donc $p'\leq 14$,  et \[
\cos\frac{\pi}{10}>\cos\frac{\pi}{11}\cos\frac{\pi}{q'} \] ce qui
impose $q'\leq 23$. On aurait alors un entier $q\geq 10$ tel que
\[\cos\frac{\pi}{10}\cos\frac{\pi}q \in
\{\cos\frac{\pi}{p'}\cos\frac{\pi}{q'}~;~11\leq p'\leq 14~;~p'\leq
q'\leq 23 \} \] ce qui est impossible.\\

\underline{\textit{Le cas $p=9$}}\\
On écarte les cas où $p'\leq 8$ par le lemme \ref{l4}. Le cas
$p'=9$ est trivial. Montrons qu'un groupe $\Gamma(2,9,q)$ ne peut
pas avoir le même spectre des longueurs qu'un groupe
$\Gamma(2,p',q')$ avec $p'\geq 10$.

L'égalité des systoles impliquerait en effet
\[
(\cos\frac{\pi}{p'})^2\leq\cos\frac{\pi}{p'}\cos\frac{\pi}{q'}=\cos\frac{\pi}9\cos\frac{\pi}q<\cos\frac{\pi}9
\]
donc $p'\leq 12$. De même, on aurait
\[\cos\frac{\pi}9>\cos\frac{\pi}{p'}\cos\frac{\pi}{q'}\geq\cos\frac{\pi}{10}\cos\frac{\pi}{q'}\]
ce qui impose $q'\leq 20$. Il s'agit donc de trouver $q\geq 9$ tel
que
\[\cos\frac{\pi}9\cos\frac{\pi}q \in
\{\cos\frac{\pi}{p'}\cos\frac{\pi}{q'}~;~10\leq p'\leq 12~;~p'\leq
q'\leq 20 \} \] ce qui s'avère impossible.\\

\underline{\textit{Le cas $p=8$}}\\
Supposons $\Gamma(2,8,q)$ et $\Gamma(2,p',q')$ isospectraux. On a
nécessairement $p'\leq 8$ et $p'\neq 3$ par le lemme \ref{l4}. Si
$p'=8$, cela impose évidemment $q'=q$. Examinons le cas
$p'\in\{4,5,6,7\}$. L'égalité des systoles se traduit par
\[\cos\frac{\pi}{8}\cos\frac{\pi}{q}=\cos\frac{\pi}{p'}\cos\frac{\pi}{q'}<\cos\frac{\pi}7
\] ce qui impose $q\leq 14$. Il s'agit donc de trouver $q'$ entier
tel que
\[\cos\frac{\pi}{8}\cos\frac{\pi}{q}=\cos\frac{\pi}{p'}\cos\frac{\pi}{q'}~;~4\leq p'\leq
7~;~8\leq q\leq 14 \] Les calculs montrent qu'un tel $q'$ n'existe
pas.\\

\underline{\textit{Le cas $p=7$}}\\
Supposons $\Gamma(2,7,q)$ et $\Gamma(2,p',q')$ isospectraux. On a
nécessairement $p'\leq 8$ par le lemme \ref{l4}. Le cas $p'=8$
s'écarte par symétrie (cf. le cas $p=8$). Si $p'=7$, cela impose
$q'=q$. Examinons le cas $p'\in\{4,5,6\}$. L'égalité des systoles
se traduit par
\[\cos\frac{\pi}{7}\cos\frac{\pi}{q}=\cos\frac{\pi}{p'}\cos\frac{\pi}{q'}<\cos\frac{\pi}6
\] ce qui impose $q\leq 11$. Il s'agit donc de trouver $q'$ entier
tel que
\[\cos\frac{\pi}{7}\cos\frac{\pi}{q}=\cos\frac{\pi}{p'}\cos\frac{\pi}{q'}~;~4\leq p'\leq
6~;~7\leq q\leq 11 \] Là encore, un tel $q'$ n'existe pas.\\

\underline{\textit{Le cas $p=6$}}\\
Supposons $\Gamma(2,6,q)$ et $\Gamma(2,p',q')$ isospectraux. On a
là encore $p'\leq 8$ en vertu du lemme \ref{l4}. Les cas où
$p'\in\{7,8\}$ s'excluent par symétrie (cf. les cas où
$p\in\{7,8\}$) et le cas $p'=6$ impose $q'=q$. Examinons les cas
où $p'\in\{4,5\}$. L'égalité des systoles se traduit toujours par
\[\cos\frac{\pi}{6}\cos\frac{\pi}{q}=\cos\frac{\pi}{p'}\cos\frac{\pi}{q'}<\cos\frac{\pi}5
\] ce qui impose $q\leq 8$. Il s'agit donc de trouver $q'$ entier
tel que
\[\cos\frac{\pi}{6}\cos\frac{\pi}{q}=\cos\frac{\pi}{p'}\cos\frac{\pi}{q'}~;~4\leq p'\leq
5~;~6\leq q\leq 8 \] Les calculs montrent l'inexistence d'un tel $q'$.\\

\underline{\textit{Le cas $p=5$}}\\
Supposons $\Gamma(2,5,q)$ et $\Gamma(2,p',q')$ isospectraux. On a
nécessairement $p'\leq 8$ par application du lemme \ref{l4}. On
exclut les cas où $p'\in\{6,7,8\}$ par symétrie. Si $p'=5$, cela
impose $q'=q$. Examinons le cas $p'=4$. L'égalité des systoles
s'écrit
\[\cos\frac{\pi}{5}\cos\frac{\pi}{q}=\cos\frac{\pi}{4}\cos\frac{\pi}{q'}<\cos\frac{\pi}4
\] ce qui impose $q\leq 6$. Il s'agit de trouver $q'$ entier
tel que
\[\cos\frac{\pi}{5}\cos\frac{\pi}{q}=\cos\frac{\pi}{4}\cos\frac{\pi}{q'}~;~5\leq q\leq 6 \]
ce qui est là encore impossible. Reste le cas $p'=3$. Pour
celui-ci, on doit trouver $q$ et $q'$ tels que
\[
2\cos\frac{\pi}5\cos\frac{\pi}q=2(\cos\frac{\pi}{q'})^2-\frac12<\frac32
\] et cela impose $q\leq 8$. On doit résoudre \[
2(\cos\frac{\pi}{q'})^2-\frac12 \in
\{2\cos\frac{\pi}5\cos\frac{\pi}q~;~q=5,6,7,8 \} \] ce qui n'est
possible que pour $q'=10$ quand $q=5$.

Les groupes $\Gamma(2,5,5)$ et $\Gamma(2,3,10)$ ont donc les mêmes
systoles. Ils ne sont cependant pas isospectraux, car la seconde longueur du spectre est différente.\\

\underline{\textit{Le cas $p=4$}}\\
Supposons $\Gamma(2,4,q)$ et $\Gamma(2,p',q')$ isospectraux. Le
lemme \ref{l4} montre que l'on a nécessairement $p'\leq 8$. On
peut écarter les cas où $p'\in\{5,6,7,8\}$ par symétrie en se
référant aux calculs effectués dans les cas où $p\in\{5,6,7,8\}$
qui figurent plus bas dans ce paragraphe. Si $p'=4$, cela implique
$q'=q$. Examinons donc le cas $p'=3$, où l'égalité des systoles se
traduit par
\[ \sqrt{2}>2\cos\frac{\pi}4\cos\frac{\pi}q=2(\cos\frac{\pi}{q'})^2-\frac12
\] qui impose $q'\leq 15$. Reste à  résoudre \[
\sqrt{2}\cos\frac{\pi}{q} \in
\{2(\cos\frac{\pi}{q'})^2-\frac12~;~q'=7,8,\dots,15 \} \] ce qui
n'est possible que pour $q'=12$ quand $q=12$. Les groupes
$\Gamma(2,4,12)$ et $\Gamma(2,3,12)$ ont donc les mêmes systoles.
Ils ne sont cependant pas isospectraux, ce qui se constate là
encore en considérant la seconde longueur du spectre.\\

\underline{\textit{Le cas $p=3$}}\\
Le lemme \ref{l4} montre que $\Gamma(2,3,q)$ ne peut pas être
isospectral à $\Gamma(2,p',q')$ avec $p'\geq 6$. Les cas $p'=4$ et
$p'=5$ s'avèrent impossibles également (cf. les calculs effectués
pour les cas $p=4$ et $p=5$). Enfin,
$p'=3$ impose clairement $q'=q$.\\

\underline{\textit{Le cas $p\geq 11$}}\\
Le cas $p'\geq 11$ est exclu par la proposition \ref{c_2}. Les cas
où $p'\leq 8$ sont impossibles cette fois ci en invoquant le lemme
\ref{l4}. Reste les cas où $p'\in\{9,10\}$ qui doivent être écartés
par symétrie (on renvoie aux calculs effectués dans les cas
$p\in\{9,10\}$).

Ceci achève la preuve du théorème \ref{t_2}. $\square$

\end{document}